\documentclass[12pt]{amsart}
\usepackage{a4wide} 
\usepackage{amsmath} 
\usepackage{amsfonts}
\usepackage{amsthm} 
\usepackage{amssymb} 
\usepackage[mathscr]{euscript}
\usepackage[all]{xy} 
\theoremstyle{plain}
\newtheorem{thm}{Theorem}
 
\newtheorem{lem}[thm]{Lemma} 
\newtheorem{cor}[thm]{Corollary}

\theoremstyle{remark}

\newtheorem*{example*}{Example}

\newtheorem*{rem*}{Remark}

\newcommand{\N}{\mathbb{N}}
\newcommand{\Z}{\mathbb{Z}}
\newcommand{\Q}{\mathbb{Q}}
\newcommand{\R}{\mathbb{R}}
\newcommand{\C}{\mathbb{C}}

\renewcommand{\H}{\mathbb{H}}

\newcommand{\E}{\mathbb{E}}

\newcommand{\calC}{\mathcal{C}}

\newcommand{\calH}{\mathcal{H}}

\newcommand{\calK}{\mathcal{K}}

\newcommand{\calM}{\mathcal{M}}
\newcommand{\calN}{\mathcal{N}}

\newcommand{\eps}{\varepsilon}
\newcommand{\bs}{\backslash}

\newcommand{\ord}{\operatorname{ord}}
\newcommand{\vol}{\operatorname{vol}}

\newcommand{\Log}{\operatorname{Log}}

\newcommand{\Sl}{\operatorname{SL}}
\newcommand{\PSl}{\operatorname{PSL}}
\newcommand{\Gl}{\operatorname{GL}}

\newcommand{\Mp}{\operatorname{Mp}}
\newcommand{\Orth}{\operatorname{O}}
\newcommand{\Unit}{\operatorname{U}}

\newcommand{\supp}{\operatorname{supp}}

\renewcommand{\div}{\operatorname{div}}

\begin{document}


\title[Two applications of the curve lemma]{
Two applications of the curve lemma for orthogonal groups}

\author[Jan H.~Bruinier]{Jan Hendrik Bruinier}
\date{January 10, 2003}
\address{Mathematisches Institut, Universtit\"at Heidelberg, Im Neuenheimer Feld 288, D-69120 Heidelberg, Germany}
\email{bruinier@mathi.uni-heidelberg.de} 
\thanks{The author is supported by a Heisenberg-Stipendium of the DFG}
\maketitle

\section{Introduction}

Recall that if $f$ is a meromorphic  elliptic modular form of weight $k$ for the group $\Sl_2(\Z)$, then for instance as a consequence of the Riemann-Roch theorem on modular curves we have 
\begin{align} \label{k/12}
\sum_{\tau\in \Sl_2(\Z)\bs \H} 
\frac{\ord_\tau(f)}{w_\tau} +\ord_\infty(f)= \frac{k}{12}.
\end{align}
Here $\H$ denotes the complex upper half plane and $w_\tau$ the order of the stabilizer of $\tau$ in $\PSl_2(\Z)$.
One purpose of this note is to give a generalization of this formula to modular forms on the orthogonal group $\Orth(2,p)$, which only depends on the Baily-Borel compactification of the arithmetic quotient in question.
Moreover, we show that certain integrals, occurring in Arakelov intersection theory, associated with modular forms on $\Orth(2,p)$ converge.
 
\medskip

We now describe the content in more detail.
Let $(V,Q)$ be a non-degenerate real quadratic space of signature $(2,p)$ and write $G\cong \Orth(2,p)$ for its real orthogonal group. 
We denote the bilinear form corresponding to the quadratic form $Q(\cdot)$ by $(\cdot,\cdot)$ such that $Q(x)=\frac{1}{2}(x,x)$. Recall that the Hermitean symmetric space corresponding to $G$ can be realized as follows: Let $V_\C=V\otimes_\R\C$ be the complexification of $V$ and $P(V_\C)$ the associated projective space. We extend the bilinear form $(\cdot,\cdot)$ $\C$-bilinearly to $V_\C$ and put
\[
\calK=\{ [Z]\in P(V_\C);\quad (Z,Z)=0,\quad (Z,\bar Z)>0\}.
\]
The set $\calK$ has two connected components. We chose one of them and denote it by $\calH$. The connected component $G^0$ of the identity of $G$ acts transitively on $\calH$. The stabilizer $K$ of a base point is a maximal compact subgroup of $G^0$ and the Hermitean symmetric space  $G^0/K$ is isomorphic to $\calH$.

Let $L\subset V$ be an even lattice, $L'$ its dual, and write $\Orth(L)$ for the integral orthogonal group of $L$. 
Throughout we let $\Gamma\subset \Orth(L)\cap G^0$ be a subgroup of finite index that acts freely on $\calH$.
By the theory of Baily-Borel \cite{BB}, the quotient $Y_\Gamma=\Gamma\bs \calH$ is a (non-singular) quasi-projective algebraic variety over $\C$ of dimension $p$. It can be compactified by adding finitely many curves and points. The resulting normal complex space is denoted by $X_\Gamma$.

Let $\widetilde{\calH}\subset V_\C-\{0\}$ be the cone over $\calH$.
Let $k\in \Z$ and $\chi$ a character of $\Gamma$.  
A meromorphic 
function $F$ on $\widetilde{\calH}$ is called a meromorphic  {\em modular form} of weight $k$ and character $\chi$ for the group $\Gamma$, if
\begin{enumerate} 
\renewcommand{\labelenumi}{\roman{enumi})}      
\item $F$ is homogeneous of degree $-k$, i.e.~$F(c Z)= c^{-k}F(Z)$ for any $c\in \C-\{0\}$;
\item $F$ is invariant under $\Gamma$, i.e.~$F(g Z)=\chi(g) F(Z)$ for any $g\in \Gamma$;
\item $F$ is meromorphic  at the boundary.
\end{enumerate}

By the Koecher principle the third condition is automatically fulfilled if the Witt rank of $L$, i.e.~the dimension of a maximal isotropic subspace of $L\otimes_\Z\Q$, is smaller than $p$. (Note that because of the signature the Witt rank of $L$ is alway $\leq 2$.)

Meromorphic  modular forms of weight $k$ with trivial character can be viewed as global rational sections of an algebraic line bundle $\calM_k$ on $Y_\Gamma$.
If $F$ is a modular form of weight $k$, then its {\em Petersson metric} is the function on $\calH$ given by 
\[
\| F(Z)\|^2_{Pet}= | F(Z) |^2 (Z,\bar Z)^k.
\]
This defines a Hermitean metric on the line bundle $\calM_k$. 

Let $\Omega$ be the first Chern form of $\calM_1$.
It is the $(1,1)$-form corresponding to the (up to a positive multiple) unique $G^0$-invariant K\"ahler metric on $\calH$.
We denote by $\vol(Y_\Gamma)=\int_{Y_\Gamma}\Omega^p$ 
the volume of $Y_\Gamma$.
If $D$ is a divisor on $Y_\Gamma$, then we define its {\em analytic degree} by
\[
\deg_{Y_\Gamma}(D)=\int_{D} \Omega^{p-1}.
\]
If $D$ is an effective divisor, then the positivity of $\Omega$ implies that $\deg_{Y_\Gamma}(D)\in \R_{\geq 0}\cup \infty$. 
Moreover, $\deg_{Y_\Gamma}(D)$ vanishes in this case, if and only if $D=0$.
(The proof of Lemma \ref{boundary3} will actually show that $\deg_{Y_\Gamma}(D)$ is finite for any divisor $D$ on $Y_\Gamma$.)

In section \ref{sect4} we prove the following theorem, which can be viewed as a higher dimension analogue of \eqref{k/12}.

\begin{thm}\label{intro2}
Assume that  the Witt rank of $L$ is $0$ if $p=1$, and $\leq 1$ if $p=2,3$ (and no restriction otherwise). If $F$ is a meromorphic modular form of weight $k$ with 
some character for the group $\Gamma$, then 
\[
\deg_{Y_\Gamma}(\div(F)) = k \vol(Y_\Gamma).
\]
\end{thm}

If $Y_\Gamma$ is compact (which can only occur if $p\leq 2$), this immediately follows from the Poincar\'e-Lelong formula. However, if $Y_\Gamma$ is non-compact, it is not obvious at all, because of the singularity of $X_\Gamma$ and the metric $\|\cdot\|_{Pet}$ at the boundary.
For modular forms whose divisor is a linear combination of Heegner divisors, 
the above formula also follows from \cite{Ku} (2.31) and \cite{Br1} Corollary 4.24.
Notice that for $p\leq 2$ the hypothesis on the Witt rank of $L$ is clearly needed. 
For instance, we may view the classical Delta function $\Delta$ as a modular form on $\Orth(2,1)$ with respect to the orthogonal group of a lattice with Witt rank $1$. By \eqref{k/12} we also have to take account of the zero of $\Delta$ at the cusp $\infty$.
However, for $p=3$ we do not have any counter-example, and it would be interesting to see whether the hypothesis on the Witt rank can be dropped. 

If $H$ is a Heegner divisor on $Y_\Gamma$, then $\deg_{Y_\Gamma}(H)$ can be computed explicitly (see \cite{Ku}, \cite{BrKue}). It is given by a Fourier coefficient of a certain Eisenstein series of weight $1+p/2$ for the metaplectic group $\Mp_2(\R)$, and can be evaluated in terms of special values of Dirichlet $L$-series and generalized divisor sums. Therefore Theorem \ref{intro2} has useful applications in this context.
%

In section \ref{sect3} we show:

\begin{thm}\label{intro1}
Assume that  the Witt rank of $L$ is $0$ if $p=1$, and $\leq 1$ if $p=2$ (and no restriction otherwise).
If $F$ is a meromorphic modular form with some character 
for the group $\Gamma$, then 
$\log \| F\|_{Pet}$
is in $L^1(Y_\Gamma,\Omega^p)$, i.e.~the integral
$
\int_{Y_\Gamma} \big| \log \| F\|_{Pet}\big| \,\Omega^p
$
converges.
\end{thm}

In the same way as with Theorem \ref{intro2}, the assumption on the Witt rank of $L$ is crucial, as for instance the example of the Delta function shows. 
Integrals of this type naturally occur in the arithmetic intersection theory of (regular models) of Shimura varieties (see e.g.~\cite{Ku} and \cite{BBK}).

The proof of the two theorems relies on the ``curve lemma'' for orthogonal groups (Theorem \ref{cl} here, see \cite{Fr1} Satz 5.8 and \cite{Fr2} Satz 2.1 for the similar case of Siegel modular varieties).
It can be used to approximate holomorphic charts of any desingularization of $X_\Gamma$ at the Baily-Borel boundary.

It seems plausible that the results of this note hold in greater generality for arithmetic quotients of Hermitean symmetric domains (under some condition on the codimension of the Baily-Borel boundary). However, Heegner divisors (see \cite{KM}, \cite{Lo}, \cite{Br1}) only exist for such domains associated with $\Orth(2,p)$ or $\Unit(1,n)$, which is why we restricted ourselves to the $\Orth(2,p)$ case.

I thank E.~Freitag for explaining the curve lemma to me and for various useful discussions on this note. Moreover, I thank U.~K\"uhn for his help.

\section{Preliminaries}
\label{sect2}

If $R>0$, we let $\E(R)=\{z\in \C;\; |z|<R\}$. For the unit disc $\E(1)$ we simply write $\E$ and put $\dot \E=\E-\{0\}$. The upper complex half plane is denoted by $\H=\{z\in\C;\; \Im(z)>0\}$. Moreover, if $f$ and $g$ are complex numbers or complex valued functions on some domain, then we write $f\ll g$, if there exists a constant $C>0$ such that $|f|\leq C |g|$.

Sometimes it is very useful to work with a tube domain realization of the Hermitean symmetric space $\calH$.
Let $z\in L$ be a primitive isotropic vector. Then we may chose $z'\in L'$ such that $(z,z')=1$. The lattice $K=L\cap z^\perp \cap {z'}^\perp$ is Lorentzian, i.e.~has signature $(1,p-1)$. 
Let $\calC$ be the cone of positive norm vectors in $K\otimes_\Z\R$.
Then $\calK_z=\{Z\in K\otimes_\Z\C;\; \Im(Z)\in \calC\}$ is isomorphic to $\calK$ via the mapping $Z\mapsto [Z+z'-Q(Z)z-Q(z')z]$.
The cone $\calC$ consists of two connected components. One of these, denoted by $\calC^+$, has the property that
\[
\calH_z= \{Z\in K\otimes_\Z\C;\; \Im(Z)\in \calC^+\}
\]
is mapped isomorphically to $\calH$ via the above mapping.
The connected component of the orthogonal group of $K\otimes_\Z\R$ viewed as a subgroup of $G^0$ acts linearly on $\calH_z$.
If $Z\in \calH_z$, then we will often denote $X=\Re(Z)$ and $Y=\Im(Z)$.
A modular form $F$ of weight $k$ for the group $\Gamma$, can be viewed as a function on $\calH_z$ with a certain transformation behavior under $\Gamma$ (see e.g.~\cite{Br1}). The Petersson metric of $F$ can be written as $\|F(Z)\|_{Pet}^2=| F(Z) |^2 (4Q(Y))^k$.

On $\calH_z$ the K\"ahler form $\Omega$ is given by
\[
\Omega=-dd^c \log (Q(Y))
\]
with $d^c=\frac{1}{4\pi i}(\partial-\bar\partial)$.
Let $e_1,\dots,e_p$ be a basis of $K\otimes_\Z\R$. If $X\in K\otimes_\Z \R$, then we
write $X=(x_1,\dots,x_p)=x_1 e_1+\dots x_p e_p$ with $x_i\in \R$.
It is easily checked that in these coordinates $d Q(Y) = (Y,dY)$ and $d^c Q(Y)=-\frac{1}{4\pi}(Y,dX)$, where $dX$ in the scalar product stands for the vector $dX=(dx_1,\dots,dx_p)$, and $dY$ for the analogous vector in the $Y$-coordinates. 
It follows that
\[
\Omega=\frac{1}{4\pi Q(Y)}\left( \frac{(Y,dX)(Y,dY)}{Q(Y)}-(dX,dY)\right).
\]
The $(p-1,p-1)$-form  $\Omega^{p-1}$ on $\calH_z$ is also easily computed explicitly.
If we put 
\begin{align}\label{dhat}
\widehat{dx_j}&=(-1)^j dx_1\cdots dx_{j-1} dx_{j+1} \cdots dx_p,
\end{align}
and define $\widehat{dy_j}$ analogously, then
\begin{align}\label{omegap-1}
\Omega^{p-1} = \frac{C}{Q(Y)^p}\sum_{i,j=1}^p \left( y_i y_j - s^{ij}Q(Y)\right)\widehat{dx_i}\widehat{dy_j},
\end{align}
where $S=(s_{ij})$ denotes the Gram matrix of $Q$ with respect to the basis $e_1,\dots,e_p$, $S^{-1}=(s^{ij})$ its inverse, and $C=(\frac{-1}{2\pi})^{p-1}\frac{(p-1)!}{2}$. 

The volume form $\Omega^p$ is up to a constant multiple equal to $\frac{DX\,DY}{Q(Y)^p}$, where $DX=dx_1\cdots dx_p$ and $DY=dy_1\cdots dy_p$ (see also \cite{BrKue}).

We will need the following elementary lemmas on vectors in the closure $\overline \calC^+\subset K\otimes_\Z\R$ of the positive cone $\calC^+$.

\begin{lem}\label{cone}
a) If $\lambda_1,\lambda_2\in \overline\calC^+$, then $(\lambda_1,\lambda_2)\geq 0$ and
\[
Q(\lambda_1)+Q(\lambda_2)\leq Q(\lambda_1+\lambda_2)\leq 2Q(\lambda_1)+2Q(\lambda_2).
\]
b) If $\lambda_1,\lambda_2\in \overline\calC^+$ are linearly independent (hence both non-zero) and $\lambda_1$ is isotropic, then  $(\lambda_1,\lambda_2)> 0$.
\end{lem}

The proof is left to the reader.

\begin{lem}\label{cone2}
Let $\lambda_1,\dots,\lambda_n\in \overline\calC^+$ be linearly independent vectors.


a) If $n\geq 2$, then for any 
$0\leq A \leq 1-2/n$ we have 
\begin{equation}\label{cone21}
Q(t_1 \lambda_1 +\dots+t_n \lambda_n)\gg t_1^{2/n+A} (t_2 \cdots t_n)^{2/n-A/(n-1)},
\end{equation}
uniformly for $t_1,\dots,t_n\in \R_{\geq 0}$.

b) Moreover, if $Q(\lambda_1)>0$, then \eqref{cone21} actually holds for any $0\leq A \leq 2-2/n$ (and trivially also for $n=1$ with $A=0$). 
\end{lem}

This can be proved by induction on $n$ using Lemma \ref{cone} and the inequality between arithmetic and geometric mean.

\section{The curve lemma}
\label{sect2.5}

We begin by recalling some facts on the Baily-Borel compactification  of $Y_\Gamma$ (see \cite{Lo} and \cite{BF} for more details).
The zero-dimensional boundary components of $Y_\Gamma$ correspond to $\Gamma$-classes of one-dimensional isotropic subspaces of $L\otimes\Q$.
The one-dimensional boundary components of $Y_\Gamma$ correspond to $\Gamma$-classes of two-dimensional isotropic subspaces of $L\otimes\Q$.

If $B$ is a one-dimensional isotropic subspace of $L\otimes\Q$, then
$B\otimes\C$ defines a point in the zero-quadric $\calN=\{ [Z]\in P(V_\C);\ (Z,Z)=0\}\subset P(V_\C)$. This is the zero-dimensional rational boundary component of $\calH$ corresponding to $B$.
If $B$ is a two-dimensional isotropic subspace of $L\otimes\Q$, then
$[B\otimes\C]\subset  \calN$.
A subset of $[B\otimes\C]$, isomorphic to $\H$, belongs to the closure of $\calH$ in $P(V_\C)$. This is the one-dimensional rational boundary component of $\calH$ corresponding to $B$.
We write $\calH^*$ for the union of $\calH$ with all zero- and one-dimensional rational boundary components of $\calH$. The group $\Gamma$ clearly acts on $\calH^*$.

By the theory of Baily-Borel there is a certain topology on $\calH^*$
 such that the quotient $X_\Gamma=\Gamma\bs\calH^*$, together with the quotient topology, has a natural structure as a projective algebraic variety, containing $Y_\Gamma$ as a Zariski-open subset.

If $s\in \calH^*$ is any point and $\Gamma_s$ its stabilizer in $\Gamma$, 
then the natural map   
\begin{align}\label{boundary}
\Gamma_s\bs\calH^*\longrightarrow \Gamma\bs\calH^*
\end{align}
defines an open embedding of a small neighborhood of the image of $s$.
If the group $\Gamma$ is neat\footnote{
Recall that a subgroup $\Gamma'$ of an algebraic group $H\subset \Gl(n,\C)$ is neat, if for every $g\in \Gamma'$ the multiplicative subgroup of $\C^\times$ generated by the eigenvalues of $g$ has no torsion. 
Every arithmetic subgroup of $H$ contains a congruence subgroup of finite index, which is neat (see \cite{Bo} Proposition 17.4).},
%
and $s$ is contained in the rational one- or two-dimensional boundary component of $\calH$ corresponding to an isotropic subspace $B\subset L\otimes\Q$, then $\Gamma_s$ centralizes $B$, i.e.~acts as the identity on it. In particular, if $z\in L\cap B$ is a primitive isotropic vector, then $\Gamma_s$ is contained in the stabilizer $\Gamma_z$ of $z$ in $\Gamma$.


The group $\Gamma_z$ can be described as follows.
We chose a vector $z'\in L'$ with $(z,z')=1$ and write $K$ for the lattice $L\cap z^\perp\cap {z'}^\perp$ as in section \ref{sect2}. 
We let $\Delta(L)\subset \Orth(L)$ be the {\em discriminant kernel} of $\Orth(L)$, i.e., the kernel of the natural homomorphism $\Orth(L)\to \Orth(L'/L)$. By our assumption, $\Gamma$ is  commensurable with $\Delta(L)$. Hence $\Gamma_z$ is commensurable with the stabilizer $\Delta(L)_z$, and it suffices (for our purposes) to consider the latter group.
We first notice that the natural inclusion of $\Orth(K\otimes \Q)$ into $\Orth(L\otimes \Q)$ induces an inclusion of the discriminant kernels $\Delta(K)\to\Delta(L)$ (but in general not of the integral orthogonal groups).
Clearly $\Delta(K)$ stabilizes $z$.
A second type of elements is given by translations. For any $\lambda\in L\cap z^\perp$, the Eichler transformation
\[
x\mapsto E(z,\lambda)(x)=x+(x,z)\lambda-(x,\lambda)z-Q(\lambda)(x,z)z \qquad(x\in V)
\] 
belongs to $\Delta(L)_z\cap G^0$ and is unipotent.
We identify $K$ with its corresponding group of Eichler transformations. The group $\Delta(K)$ acts on $K$ and we claim that $\Delta(L)_z=\Delta(K) \ltimes K$.
In fact, let $g\in \Delta(L)_z$. Since $g$ fixes $L'/L$, the vector $z'-gz'$ is contained in $L$. It is clearly orthogonal to $z$ and therefore 
$E(z,z'-gz')\in \Delta(L)_z$. One easily checks that
$E(z,z'-gz') g$ stabilizes $z'$ and $z$, and consequently is contained in 
$\Delta(L)_z\cap \Orth(K\otimes \Q)=\Delta(K)$. This proves the claim.
The action of $g=(U,T)\in\Delta(K) \ltimes K$ on the tube domain $\calH_z$ is given by
\[
g Z = U Z + T \qquad(Z\in \calH_z).
\]

The following result is the fundamental for our later argument.

\begin{thm}\label{cl}
Let $\Psi: \H^m\times \E^N\to \calH_z$ be a holomorphic mapping satisfying
\begin{align*}
\Psi(\tau_1+1,\tau_2,\dots,\tau_m,w)&=U_1 \Psi(\tau_1,\dots,\tau_m,w) + T_1,\\
 &\;\;\vdots\\
\Psi(\tau_1,\tau_2,\dots,\tau_m+1,w)&=U_m \Psi(\tau_1,\dots,\tau_m,w) + T_m,
\end{align*}
where $U_1,\dots,U_m\in \Delta(K)$, and $T_1,\dots,T_m$ are translations in $K$, and $w\in \E^N$.
Then  $U_1,\dots,U_m$ are of finite order $\ell$ for some $\ell\in \N$, and there exist semi-positive vectors $\lambda_1,\dots,\lambda_m\in (K\otimes\Q) \cap \overline\calC^+$ and a holomorphic function $\Psi_0: \E^m\times \E^N\to \calH_z$  such that
\[
\Psi(\tau_1,\dots,\tau_m,w)=\lambda_1\tau_1+\dots+\lambda_m \tau_m +\Psi_0(e^{2\pi i \tau_1/\ell},\dots,e^{2\pi i \tau_m/\ell},\, w).
\]
\end{thm}

\begin{proof}
Theorem \ref{cl} for $m=1$ is the orthogonal group analogue of the curve Lemma for Siegel modular groups (see \cite{Fr1} Satz 5.8, \cite{Fr2} Satz 2.1).
It can be proved in a similar way. The case $m>1$ can be reduced to the previous case by considering the functions $(\tau_i,w) \mapsto \Psi(\tau_1,\dots,\tau_m,w)$ for $i=1,\dots,m$, where the $\tau_j$-variables with $j\neq i$ are being fixed.
\end{proof}

\section{The integral of the logarithm of a modular form}
\label{sect3}




\begin{proof}[Proof of Theorem \ref{intro1}]
Without any restriction we may assume that $\Gamma$ is neat and acts trivially on $L'/L$. 

Since $X_\Gamma$ is compact, it suffices to show that $\log \| F\|_{Pet}$ is locally integrable near every point of $X_\Gamma$.
This is easily done for the points of $Y_\Gamma$. Thus we only have to consider the points of the Baily-Borel boundary $\partial X_\Gamma = X_\Gamma-Y_\Gamma$. In particular, in the case that $Y_\Gamma$ is already compact we are done.
Because of our assumption on the Witt rank of $L$ we may therefore assume that $p\geq 2$.

According to Hironaka's theory there exists a desingularization $\pi:\widetilde X_\Gamma\to X_\Gamma$ of $X_\Gamma$ with respect to the divisor $\div(F)$ such that $\pi^{-1}(\div(F)\cup \partial X_\Gamma)$ is a divisor with normal crossings.
Let $a\in \partial X_\Gamma$ be a boundary point and $\tilde a\in \widetilde X_\Gamma$ be a point with $\pi(\tilde a)=a$. 
Let $U\subset \widetilde X_\Gamma$ and $V\subset X_\Gamma$ be open neighborhoods of $\tilde a$ and $a$, respectively, such that $\pi$ induces an isomorphism
$U-\pi^{-1}(\div(F)\cup \partial X_\Gamma)\to V-(\div(F)\cup \partial X_\Gamma)$.
After a biholomorphic change of coordinates we may assume without loss of generality that $U=\E^p$, $\tilde a=(0,\dots ,0)\in \E^p$, and 
\begin{align}\label{divboundary}
U\cap\pi^{-1}(\partial X_\Gamma)&=\{(q_1,\dots,q_p)\in \E^p;\; q_1\cdots q_m=0\}\\
\label{divf}
U\cap\pi^{-1}(\div(F))&=\{(q_1,\dots,q_p)\in \E^p;\; q_{i_1}\cdots q_{i_u}=0\}
\end{align}
for some $1\leq m\leq p$ and $1\leq i_1<\dots<i_u\leq p$. Moreover, we may assume that $V$ is given by the open embedding of a neighborhood of a boundary point $s\in \calH^*$ above $a$ as in \eqref{boundary}.

We obtain a holomorphic mapping of complex spaces
\[
\psi: \dot \E^m\times \E^{p-m} \longrightarrow \Gamma_s\bs\calH,
\]
which can be continued to a holomorphic mapping $ \E^{p} \to \Gamma_s\bs\calH^*$ whose image is $V$.
We lift $\psi$ to a holomorphic mapping $\Psi: \H^m\times \E^{p-m} \longrightarrow \calH$
of the universal covers and get the following commutative diagram:
\begin{align}\label{diag}
\xymatrix{  \H^m\times\E^{p-m} \ar[d]  \ar[r]^{\Psi} &    \calH \ar[d]\\
\dot{\E}^m\times \E^{p-m}\ar[d]   \ar[r]^{\psi}  &\Gamma_s\bs\calH\ar[d]  \\
\E^m\times \E^{p-m}  \ar[r] &\Gamma_s\bs\calH^*}.
\end{align}
The local integrability of  $\log\|F\|_{Pet}$ near $a\in X_\Gamma$  is equivalent to showing that
\begin{align}\label{int2}
\int\limits_{\E(R)^p} \big| \psi^*(\log\|F\|_{Pet})\big| \,\psi^*(\Omega^p)
\end{align}
converges for some $0<R<1$.

Throughout we use $\tau=(\tau_1,\dots,\tau_m)$ as a standard variable in $\H^m$,   $q=(q_1,\dots,q_p)$  as a standard variable in $\E^p$, and frequently denote $w=(q_{m+1},\dots,q_p)\in \E^{p-m}$.
Since the diagram commutes, there exist elements $g_1,\dots,g_m\in \Gamma_s$ such that
\begin{align*}
\Psi(\tau_1+1,\tau_2,\dots,\tau_m,w)&=g_1 \Psi(\tau_1,\dots,\tau_m,w),\\
 &\;\;\vdots\\
\Psi(\tau_1,\tau_2,\dots,\tau_m+1,w)&=g_m \Psi(\tau_1,\dots,\tau_m,w).
\end{align*}
The transformations $g_1,\dots,g_m$ do not depend on $(\tau,w)$, because $\Gamma$ acts properly discontinuously and freely.
Moreover, they have to fix the boundary point $s$.

Let $z\in L$ be a primitive isotropic vector such that $\Gamma_s\subset \Gamma_z$. From now on we work with the tube domain realization $\calH_z$ of $\calH$.
By virtue of the remarks of section \ref{sect2.5} we may conclude 
that any $g_i$ has the form
$g_i Z = U_i Z + T_i$,
where $U_i\in \Delta(K)$, $T_i\in K$, and $Z\in \calH_z$.
  
In view of Theorem \ref{cl} there exist semi-positive vectors $\lambda_1,\dots,\lambda_m\in (K\otimes\Q) \cap \overline\calC^+$ and a holomorphic function $\Psi_0: \E^p\to \calH_z$  such that
\[
\Psi(\tau_1,\dots,\tau_m,w)=\lambda_1\tau_1+\dots+\lambda_m \tau_m +\Psi_0(e^{2\pi i \tau_1},\dots,e^{2\pi i \tau_m},w).
\]
Here we have used our assumption that $\Gamma$ be neat and the fact that the translations $T_i$ are unipotent to conclude that $\ell=1$ and $U_i=1$.  
Consequently $\psi$ is given by 
\begin{align}\label{psi}
\psi(q)&=\frac{1}{2\pi i} \left(\lambda_1 \Log q_1 +\dots +\lambda_m \Log q_m\right) + \Psi_0(q),
\end{align}
where $\Log(\cdot)$ denotes a fixed branch of the holomorphic logarithm on $\E-\R_{\leq 0}$. 
The fact that $\psi(q)\in \partial X_\Gamma$, if $q_i=0$ for some $1\leq i\leq m$, implies that $\lambda_1,\dots,\lambda_m\neq 0$.
Using polar coordinates
\begin{align*}
q_j= r_j e^{i\rho_j} \qquad (0\leq r_j <1,\; 0\leq \rho_j < 2\pi),
\end{align*}
we obtain for the imaginary part of $\psi$:
\begin{align}\label{impsi}
\psi^*(Y)=\Im(\psi(q))&=-\frac{1}{2\pi} \left(\lambda_1 \log r_1 +\dots +\lambda_m \log r_m\right) + \Im(\Psi_0(q)).
\end{align}
In order to prove the convergence of \eqref{int2}, we have to find bounds for $\psi^*(Q(Y))$ and for the functional determinant of $\psi$. Since $\Psi_0(q)$ is bounded on $\E(R)^p$, this can be done by means of \eqref{psi} and \eqref{impsi}.

Let $\mu$ be the rank of $\lambda_1,\dots,\lambda_m$ ($1\leq \mu\leq m$). Without any restriction we may assume that $\lambda_1,\dots,\lambda_\mu$ are linearly independent.
Using \eqref{impsi} and Lemma \ref{cone}, we see that
\begin{align}\label{ul}
0< q\left(\lambda_1 \log r_1 +\dots +\lambda_\mu \log r_\mu\right) \ll \psi^*(Q(Y))\ll (\log r_1)^2 +\dots +(\log r_m)^2.
\end{align}
In view of \eqref{divf} the upper bound in particular implies that
\begin{align}\label{boundf}
|\psi^*(\log\|F\|_{Pet})|\ll |\log r_1|+\dots +|\log r_p|
\end{align}
uniformly on $\E(R)^p$.

The complex Jacobi Matrix of $\psi$ is given by
\begin{align*}
J(\psi;q)&
= \frac{1}{2\pi i}\begin{pmatrix}\lambda_1/q_1&\cdots&\lambda_m/q_m&0&\cdots&0
\end{pmatrix}
+ J(\Psi_0;q) \\
&=\Lambda\begin{pmatrix}q_1^{-1}& & &\\
& \ddots& &\\
& & q_m^{-1} & \\
& & & 0_{p-m}
\end{pmatrix} 
+J(\Psi_0;q) 
\end{align*}
where $\Lambda$ denotes the $p\times p$ matrix $\Lambda=\frac{1}{2\pi i}\begin{pmatrix}\lambda_1&\cdots&\lambda_m &0&\cdots&0
\end{pmatrix}$. The fact that $\lambda_1,\dots,\lambda_\mu$ form a basis for the span of all columns of $\Lambda$ implies that there is a matrix $S$ of elementary row manipulations such that 
\[
S\Lambda = \begin{pmatrix} 
l_{11} & * & \cdots & * \\
0   &\ddots & \ddots& \vdots\\
\vdots & \ddots &l_{\mu\mu} & *\\
0& \cdots & 0& 0
\end{pmatrix}.   
\]
The diagonal elements $l_{11},\dots,l_{\mu\mu}$ are non-zero such that the first $\mu$ rows are linearly independent. The rows $\mu+1$ to $p$ vanish identically.
We find that the determinant of $J(\psi;q)$ is up to the sign equal to the determinant of
\[
\begin{pmatrix} 
l_{11}/q_1 & * & \cdots & * \\
0   &\ddots & \ddots& \vdots\\
\vdots & \ddots &l_{\mu\mu}/q_\mu & *\\
0& \cdots & 0& 0
\end{pmatrix}
+ S J(\Psi_0;q).
\]
If we expand the determinant of this matrix successively by the first to the $\mu$-th column and use the fact that $S J(\Psi_0;q)$ is bounded, we find that
\begin{align*}
\det(J(\psi;q))&\ll r_1^{-1}\cdots r_\mu^{-1}
\end{align*}
and therefore
\begin{align}
\nonumber
\psi^*(\Omega^p)&\ll \psi^*\left(\frac{DX\, DY}{Q(Y)^p}\right)
= \psi^*(Q(Y)^{-p}) \cdot|\det J(\psi;q)|^2 \cdot 
r_1 dr_1 d\rho_1\cdots r_p dr_p d\rho_p\\
\label{boundj1}
&\ll \psi^*(Q(Y)^{-p}) 
\frac{d\rho_1\cdots d\rho_p\, dr_1\cdots dr_p}{r_1\cdots r_\mu}.
\end{align}
In view of \eqref{boundf} and \eqref{boundj1}, to prove that the integral  \eqref{int2} is finite, it suffices to show that
\begin{align}\label{intj}
\int\limits_{\E(R)^p} 
\frac{|\log r_j|}{\psi^*(Q(Y)^{p})} 
\frac{d\rho_1\cdots d\rho_p\, dr_1\cdots dr_p}{r_1\cdots r_\mu}
\end{align}
converges for $j=1,\dots,p$. The case $\mu+1\leq j\leq p$ is easily treated using \eqref{boundqy1} and left to the reader.
For the case $1\leq j\leq \mu$ we may assume without loss of generality that $j=1$. 
According to Lemma \ref{cone2}a, we may infer from the lower bound of 
\eqref{ul} (respectively directly from \eqref{impsi}, if $\mu=1$ and $\lambda_1$ is anisotropic) that 
\begin{align}\label{boundqy1}
\psi^*(Q(Y))\gg |\log r_1|\left( |\log r_2| \cdots |\log r_\mu|\right)^{1/(\mu-1)},
\end{align}
uniformly on $\E(R)^p$. Here the right hand side is understood as $|\log r_1|$, if $\mu=1$.
Hence it suffices to prove that
\begin{align}\label{h1}
\int\limits_{r_1=0}^R \cdots  \int\limits_{r_p=0}^R
\frac{ dr_1\cdots dr_p}{r_1\cdots r_\mu 
\cdot |\log r_1|^{p-1} \cdot|\log r_2|^{p/(\mu-1)} \cdots |\log r_\mu|^{p/(\mu-1)}}
\end{align}
converges.
Since
\begin{align}\label{intest1}
\int\limits_{r=0}^R r^\alpha  \frac{ dr}{|\log r|^{\beta}}<\infty, &\quad 
\text{for $\alpha>-1$ and $\beta\in \R$},\\
\label{intest2}
\int\limits_{r=0}^R \frac{ dr}{r |\log r|^\beta}<\infty, &\quad 
\text{for $\beta>1$},
\end{align}
we find that the integrals over $r_{2},\dots, r_p$ converge.
Thus \eqref{h1} is bounded by a positive constant multiple of
\[
\int\limits_{r_1=0}^R \frac{ dr_1}{r_1 |\log r_1|^{p-1}}.
\]
This is finite  if $p>2$.
However, if $p=2$, then 
our assumption on the Witt rank of $L$ implies that $K$ is anisotropic. Hence, by Lemma \ref{cone2}b, we have the estimate
\begin{align}\label{boundqy2}
\psi^*(Q(Y))\gg \begin{cases}|\log r_1|^{2},&\text{if $\mu=1$,}\\
 |\log r_1|^{1+\eps} |\log r_2|^{1-\eps},&\text{if $\mu=2$,}
                \end{cases}
\end{align}
where $\eps>0$ is small.
Consequently,
\begin{align*}
\psi^*(\Omega^p)
&\ll \begin{cases}r_1^{-1}\frac{d\rho_1  d\rho_2\, dr_1  dr_2}{|\log r_1|^{4}},&\text{if $\mu=1$,} \\
r_1^{-1} r_2^{-1} \frac{d\rho_1  d\rho_2\, dr_1  dr_2}{|\log r_1|^{2+2\eps}|\log r_2|^{2-2\eps}}, &\text{if $\mu=2$.}
     \end{cases}
\end{align*}
The same argument as above then shows that  
the integral \eqref{intj} (with $j=1$)  is actually bounded by
\[
\int\limits_{r_1=0}^R \frac{ dr_1}{r_1 |\log r_1|^{1+\delta}}
\]
for some small $\delta>0$.
This concludes the proof of the theorem.
\end{proof}

\begin{rem*} The argument of the proof of this theorem actually shows that $\log \| F\|_{Pet}$ belongs to $L^{1+\eps}(X_\Gamma,\Omega^p)$ for some $\eps>0$.
\end{rem*}

\section{The degree formula}\label{sect4}


If $Y$ is a complex manifold, we write $C^\infty(Y)$  for the space of $C^\infty$-functions on $Y$. 
If $X$ is any complex space, we define $C^\infty(X)$ to be the space of those 
functions $f:X\to \C$, for which for any holomorphic map $h:Y\to X$ from a complex manifold $Y$ to $X$ the pull back satisfies $h^*(f)\in C^\infty(Y)$.
The subspace of compactly supported functions on $X$ is denoted by $C^\infty_c(X)$. 
Using local embeddings of $X$ into some $\C^N$, one obtains many $C^\infty$-functions with compact support by pull back. They can be used to construct a partition of unity, if $X$ is paracompact. 

%
%

\begin{proof}[Proof of Theorem \ref{intro2}]
We may and will assume that $\Gamma$ is neat and acts trivially on $L'/L$.
If $Y_\Gamma$ is compact, the assertion is an immediate consequence of the Poincar\'e-Lelong formula (see \cite{GH}, \cite{SABK} Chapter II.1.4).
In view of our assumption on the Witt rank of $L$ we may therefore assume that $p\geq 2$. 

Since the Petersson metric and $\Omega^{p-1}$ are smooth on $Y_\Gamma$, 
the Poincar\'e-Lelong formula implies that for all $\sigma\in C_c^\infty(Y_\Gamma)$:
\begin{equation}\label{pl}
k\int\limits_{X_\Gamma} \Omega^p \sigma = \int\limits_{\div_{X_\Gamma}(F)} \Omega^{p-1} \sigma - \int\limits_{X_\Gamma} \log\|F\|_{Pet}^2 \Omega^{p-1} dd^c \sigma.
\end{equation}
We will show that this formula actually holds for $\sigma\in C^\infty(X_\Gamma)$. This implies the theorem by taking $\sigma=1$.

The Baily-Borel compactification $X_\Gamma$ admits partition of unity. 
Hence it suffices to prove the assertion locally, i.e., that for every $a\in X_\Gamma$ there is a neighborhood $V\subset X_\Gamma$ such that \eqref{pl} holds for all $\sigma\in C_c^\infty(V)$.
If $a\in Y_\Gamma$, this is a consequence of the Poincar\'e-Lelong formula.
We therefore assume that $a\in \partial X_\Gamma$ is a boundary point.
In this case the assertion is not at all clear, because of the singularity of the space $X_\Gamma$ and the metric $\|\cdot\|_{Pet}$.

We begin by preparing nice local charts using Theorem \ref{cl} in the same way as in the proof of Theorem \ref{intro1}. We stick to the notation introduced before. 
We chose a desingularization $\pi:\widetilde X_\Gamma\to X_\Gamma$ of $X_\Gamma$ with respect to the divisor $\div(F)$ such that $\pi^{-1}(\div(F)\cup \partial X_\Gamma)$ is a divisor with normal crossings.
Let $\tilde a\in \widetilde X_\Gamma$ be a point with $\pi(\tilde a)=a$, and let $U\subset \widetilde X_\Gamma$ and $V\subset X_\Gamma$ be open neighborhoods of $\tilde a$ and $a$, respectively, as in the proof of Theorem  \ref{intro1}.
After a biholomorphic change of coordinates we may again assume that $U=\E^p$, $\tilde a=(0,\dots ,0)\in \E^p$, and that \eqref{divboundary}, \eqref{divf} hold.
Moreover, we may assume that $V$ is given by the open embedding of a neighborhood of a boundary point $s\in \calH^*$ above $a$ as in \eqref{boundary}.
We obtain  holomorphic mappings of complex spaces
$\psi: \dot \E^m\times \E^{p-m} \to \Gamma_s\bs \calH$ and 
$\Psi: \H^m\times \E^{p-m} \to \calH$ such that the diagram \eqref{diag} commutes.
By means of the curve lemma we find that 
there exist semi-positive non-zero vectors $\lambda_1,\dots,\lambda_m\in (K\otimes\Q) \cap \overline\calC^+$ and a holomorphic function $\Psi_0: \E^p\to \calH_z$  such that
\[
\Psi(\tau_1,\dots,\tau_m,w)=\lambda_1\tau_1+\dots+\lambda_m \tau_m +\Psi_0(e^{2\pi i\tau_1},\dots,e^{2\pi i\tau_m},w).
\]
There is a unique splitting of the divisor of $\psi^*(F)$ into a sum
\[
\div(\psi^*(F))=D_1+D_2,
\]
where the support of $D_1$ is contained in $E=\{q\in \E^p;\; q_1\cdots q_m=0\}$, the exceptional divisor over $\partial X_\Gamma$, and $\supp(D_2)\cap \supp(E)\subset\{0\}$. 
We find that it suffices to show
\begin{align}\label{toshow}
k\int\limits_{\E^p} \psi^*(\Omega^p) \sigma = \int\limits_{D_2} \psi^*(\Omega^{p-1}) \sigma - \int\limits_{\E^p} \psi^*(\log\|F\|_{Pet}^2 \Omega^{p-1}) dd^c \sigma
\end{align}
for all $\sigma\in C^\infty_c(\E^p)$. We will prove this following the standard argument for the Poincar\'e-Lelong formula and making sure that the various boundary terms caused by the singularities vanish.

Let $\sigma\in C^\infty_c(\E^p)$ and write 
$\eta=\psi^*(\Omega^{p-1})\sigma$. We put $D=\{q\in \E^p;\; q_1\cdots q_p=0\}=\E^p-\dot\E^p$ and note that $\log\|F\|_{Pet}^2$ and $\eta$ are smooth outside $D$.
We have
\begin{align*}
\int\limits_{\E^p} \psi^*(\log\|F\|_{Pet}^2) dd^c \eta &
=\int\limits_{\E^p} d\left( \psi^*(\log\|F\|_{Pet}^2) d^c \eta\right) - \int\limits_{\E^p} \left(d \psi^*(\log\|F\|_{Pet}^2)\right) (d^c \eta)\\
&=\lim_{\eps \to 0} \int\limits_{\E^p-B_\eps(D)} d\left(  \psi^*(\log\|F\|_{Pet}^2)
d^c \eta\right) + \int\limits_{\E^p} \left(d^c \psi^*(\log\|F\|_{Pet}^2) \right) (d \eta)\\
&=-\lim_{\eps \to 0} \int\limits_{\partial B_\eps(D)} \psi^*(\log\|F\|_{Pet}^2)
(d^c \eta) + 
\int\limits_{\E^p} \left(d^c \psi^*(\log\|F\|_{Pet}^2) \right) (d \eta).
\end{align*}
Here $B_\eps(D)$ denotes a tubular $\eps$-neighborhood of $D$.
It follows from Lemma \ref{boundary1} below that the limit of the integral over
$\partial B_\eps(D)$ vanishes. 
Hence 
\begin{align*}
\int\limits_{\E^p} \psi^*(\log\|F\|_{Pet}^2) dd^c \eta
&=-\int\limits_{\E^p} d\left((d^c \psi^*(\log\|F\|_{Pet}^2)) \eta\right) 
+\int\limits_{\E^p} 
\left(dd^c \psi^*(\log\|F\|_{Pet}^2) \right)  \eta\\
&= \lim_{\eps \to 0} \int\limits_{\partial B_\eps(D)} (d^c \psi^*(\log\|F\|_{Pet}^2) ) \eta 
- k \int\limits_{\E^p} \psi^*(\Omega)  \eta.
\end{align*}
We will prove in Lemma \ref{boundary2} below that 
\[ 
\lim_{\eps \to 0} \int\limits_{\partial B_\eps(D)} (d^c \psi^*(\log Q(Y)) ) \eta=0.
\]
Consequently
\begin{align*}
\int\limits_{\E^p} \psi^*(\log\|F\|_{Pet}^2) dd^c \eta &
=\lim_{\eps \to 0} \int\limits_{\partial B_\eps(D)} (d^c \psi^*(\log|F|^2) ) \eta 
- k \int\limits_{\E^p} \psi^*(\Omega)  \eta.
\end{align*}
It follows from Lemma \ref{boundary1} that the first integral on the right hand side vanishes, if $\psi^*(F)$ is a an invertible holomorphic function on $\E^p$.
Thus, by linearity we may assume that $\psi^*(F)=q_\kappa$ for some $1\leq \kappa\leq p$.
It is easily checked that $d^c\log|q_\kappa|^2=\frac{1}{2\pi} d\rho_\kappa$.
In Lemma \ref{boundary3} we show that
\begin{align*}
\frac{1}{2\pi}\lim_{\eps \to 0} \int\limits_{\partial B_\eps(D)} d\rho_\kappa\wedge \eta = 
\begin{cases} 
0,&\text{if $1\leq \kappa\leq m$,}\\
\int_{q_\kappa=0} \eta,&\text{if $m< \kappa\leq p$,}
\end{cases}
\end{align*}
concluding the proof of the Theorem.
\end{proof}

\begin{lem}\label{boundary1}
For any $n=1,\dots,p$ we have
\[
\lim_{\eps \to 0} \int\limits_{\partial B_\eps(D)} |\log r_n|
(d^c \eta) =0.
\]
\end{lem}

\begin{proof}
In the same way as in \eqref{dhat} we define 
\begin{align*}
\widehat{dq_j}&=(-1)^j dq_1\cdots dq_{j-1} dq_{j+1} \cdots dq_p,
\end{align*}
and $\widehat{d\bar q_j}$ analogously.
It follows from \eqref{omegap-1} and \eqref{impsi} that the $\widehat{dq_i} \widehat{d\bar q_j}$-component of $\psi^*(\Omega^{p-1})$ is bounded by a sum over $1\leq k,l\leq p$ of $(p-1,p-1)$-forms 
of the form
\[
\frac{|\log r_1|^2+\dots +|\log r_m|^2}{\psi^*(Q(Y)^p)} \det (J_{ik} (\psi;q)) \det (\overline{J_{jl} (\psi;q)})    \widehat{dq_i} \widehat{d\bar q_j}. 
\]
Here $J_{ik} (\psi;q)$ denotes the matrix obtained from the Jacobi matrix $J(\psi;q)$ by canceling the $i$-th column and the $k$-th row.

Let $\Lambda_{ik}$ be the matrix obtained from $\Lambda=\frac{1}{2\pi i}\begin{pmatrix}\lambda_1&\cdots&\lambda_m&0&\cdots&0
\end{pmatrix}$  by canceling the $i$-th column and the $k$-th row. Chose pairwise different indices $\alpha_1,\dots,\alpha_\mu\leq m$ (with $0\leq \mu\leq m-1$ and $\alpha_t\neq i$) such that the sub-matrix of $\Lambda_{ik}$ corresponding to 
$\lambda_{\alpha_1},\dots,\lambda_{\alpha_\mu}$ consists of a basis for the vector space generated by the columns of $\Lambda_{ik}$. Then we find as in the proof of Theorem \ref{intro1} that
\[
\det (J_{ik} (\psi;q)) \ll r_{\alpha_1}^{-1} \cdots r_{\alpha_\mu}^{-1}.
\]
In the same way there are pairwise different indices $\beta_1,\dots,\beta_{\nu}\leq m$ (with $0\leq \nu\leq m-1$ and $\beta_t\neq j$) such that the sub-matrix of $\Lambda_{jl}$ corresponding to 
$\lambda_{\beta_1},\dots,\lambda_{\beta_{\nu}}$ consists of a basis for the vector space generated by the columns of $\Lambda_{jl}$, and
\[
\det (J_{jl} (\psi;q)) \ll r_{\beta_1}^{-1} \cdots r_{\beta_{\nu}}^{-1}.
\] 

Consequently, the assertion of the lemma follows if we can show that
\begin{align}\label{boundary11} 
\lim_{\eps \to 0} \int\limits_{\partial B_\eps(D)} \frac{|\log r_n|^3}{\psi^*(Q(Y)^p)} \frac{ \widehat{dq_i} \widehat{d\bar q_j}(d^c \sigma)}{ r_{\alpha_1} \cdots r_{\alpha_\mu} r_{\beta_1}\cdots r_{\beta_{\nu}}} 
  =0 
\end{align}
for all $i,j,n=1,\dots,p$. 
Since $\sigma$ has compact support, there exists a $0<R<1$ with $\supp(\sigma)\subset \E(R)^p$.
It suffices to show that the contribution to the integral coming from 
\begin{align*}
S_t(\eps,R)&=\{q\in \partial B_\eps(D)\cap \E(R)^p;\; |q_t|=\eps\}
=\{q\in   \E(R)^p;\; \text{$|q_1|,\dots, |q_p|\geq \eps$ and $|q_t|=\eps$}\}
\end{align*}
goes to zero in the limit for $1\leq t\leq p$.
Here only the $d\rho_t \widehat{dq_t}\widehat{d\bar q_t}$ component of $\widehat{dq_i} \widehat{d\bar q_j}(d^c \sigma)$ gives a non-zero contribution.
In particular the integral vanishes, if both, $i$ and $j$, are different from $t$. Without loss of generality we may therefore assume that $j= t$. Because $d^c\sigma$ is bounded, we find that it suffices to show 
\begin{align}\label{boundary12} 
\lim_{\eps \to 0} \int\limits_{S_t(\eps,R)} \frac{|\log r_n|^3 }{\psi^*(Q(Y)^p)} \frac{ d\rho_1 dr_1\cdots d\rho_{t-1} dr_{t-1}\cdot r_td\rho_t \cdot  d\rho_{t+1} dr_{t+1}\cdots  d\rho_p dr_p}{ r_{\alpha_1} \cdots r_{\alpha_\mu} } 
  =0 
\end{align}
for all $n=1,\dots,p$.

Let us first assume that $\alpha_1,\dots,\alpha_\mu$ 
are all different from $t$.
Then we may use the (poor) estimate $ \psi^*(Q(Y))\gg 1$ to infer that the integral in \eqref{boundary12} is bounded by
\begin{align}\label{hh1}
\eps \int\limits_{r_1=\eps}^R\cdots \int\limits_{r_{t-1}=\eps}^R \cdot \int\limits_{r_{t+1}=\eps}^R\cdots \int\limits_{r_p=\eps}^R 
|\log r_n|^3
\frac{dr_1\cdots dr_{t-1}\cdot dr_{t+1}\cdots dr_{p}} { r_{\alpha_1} \cdots r_{\alpha_\mu} }.
\end{align}
For any $B\in \R$ we have  
\begin{align}\label{intest3}
\int\limits_{r=\eps}^R |\log r|^B \frac{ dr}{r}= 
\begin{cases} O(|\log\eps|^{B+1}),&\text{if $B> -1$,}\\
O(\log|\log \eps|),&\text{if $B=-1$,}\\
O(1),&\text{if $B<-1$,}
\end{cases}
\end{align}
as $\eps\to 0$.
Thus \eqref{hh1} is bounded by $\eps |\log\eps|^{C}$ for some $C>0$ as $\eps\to 0$, implying \eqref{boundary12}.

%

We now assume that one of the indices $\alpha_1,\dots,\alpha_\mu$ is equal to $t$, say $\alpha_1=t$. Then in particular $t\leq m$ and without loss of generality we may assume $\alpha_1=t=1$. We only consider the case that $n=1$ as well, leaving the similar remaining case to the reader.
It follows from Lemma \ref{cone2}a with $A=1-2/\mu$ (respectively directly from \eqref{impsi} if $\mu=1$) that
\begin{align*}
\psi^*(Q(Y))&\gg Q(\lambda_{\alpha_1} \log r_{\alpha_1}  + \dots +\lambda_{\alpha_\mu} \log r_{\alpha_\mu})\\
&\gg |\log r_{\alpha_1}|\cdot \left( |\log r_{\alpha_2}|\cdots |\log r_{\alpha_{\mu}}|
\right)^{1/(\mu-1)}.
\end{align*}
Thus the integral in \eqref{boundary12} is bounded by
\begin{align*}
&\int\limits_{S_1(\eps,R)}   |\log r_{\alpha_1}|^{3-p} \left( |\log r_{\alpha_2}|\cdots |\log r_{\alpha_{\mu}}|
\right)^{-p/(\mu-1)}
\frac{   d\rho_1 \cdot d\rho_2 dr_2 \cdots d \rho_p d r_p }{ r_{\alpha_2} \cdots r_{\alpha_\mu} }\\
&\ll |\log \eps|^{3-p} \int\limits_{r_{\alpha_2}=\eps}^R\cdots \int\limits_{r_{\alpha_\mu}=\eps}^R 
\left( |\log r_{\alpha_2}|\cdots |\log r_{\alpha_{\mu}}|
\right)^{-p/(\mu-1)}
\frac{dr_{\alpha_2}\cdots dr_{\alpha_\mu}} { r_{\alpha_2} \cdots 
r_{\alpha_\mu}}\\
&\ll  |\log \eps|^{3-p}.
\end{align*}
Here, in the last line we have used \eqref{intest3}. 
Hence we obtain \eqref{boundary12}, if $p>3$. 
In the case $p\leq 3$ our assumption on the Witt rank of $L$ implies that $K$ is anisotropic. Therefore we may apply Lemma \ref{cone2}b with $A=2-2/\mu$ to infer that
$\psi^*(Q(Y))\gg |\log r_{1}|^2$.
By means of \eqref{intest3} we find that the integral in \eqref{boundary12} is actually bounded by $|\log \eps|^{2-2p+\mu}\ll|\log \eps|^{-1}$, which implies \eqref{boundary12}.
\end{proof}

\begin{lem}\label{boundary2} 
We have
\[ 
\lim_{\eps \to 0} \int\limits_{\partial B_\eps(D)} (d^c \psi^*(\log Q(Y)) ) \eta=0.
\]
\end{lem}

\begin{proof}
One can use \eqref{impsi}, the estimate $\psi^*(Q(Y))\gg |\log r_1|+\dots+|\log r_m|$, and the fact that $d^c\log|q_j|^2=\frac{1}{2\pi} d\rho_j$, to conclude that
\[
d^c \psi^*(\log Q(Y)) \ll  d\rho_1 +\dots +d\rho_m + r_{m+1}d\rho_{m+1}+\dots+r_p d\rho_p + dr_1+\dots +dr_p .
\]
Here ``$\ll$'' is understood componentwise.
Thus the assertion follows from the estimate for \eqref{boundary12} in Lemma \ref{boundary1}.
\end{proof}

\begin{lem}\label{boundary3} 
We have
\begin{align*}
\frac{1}{2\pi}\lim_{\eps \to 0} \int\limits_{\partial B_\eps(D)} d\rho_\kappa\wedge \eta = 
\begin{cases} 
0,&\text{if $1\leq \kappa\leq m$,}\\
\int_{q_\kappa=0} \eta,&\text{if $m< \kappa\leq p$.}
\end{cases}
\end{align*}
In the latter case, the integral on the right hand side converges absolutely.
\end{lem}

\begin{proof}
We first consider the case that $1\leq \kappa\leq m$, say $\kappa=1$.
By the argument of Lemma \ref{boundary1} (using the notation of that lemma) it suffices to show that
\begin{align}\label{boundary31}
\lim_{\eps \to 0} \int\limits_{S_t(\eps,R)} \frac{|\log r_n|^2 }{\psi^*(Q(Y)^p)} \frac{  d\rho_1  \widehat{dq_i} \widehat{d\bar q_j}}{ r_{\alpha_1} \cdots r_{\alpha_\mu} r_{\beta_1}\cdots r_{\beta_{\nu}}} 
  =0 
\end{align}
for all $i,j,t=1,\dots,p$ and $n=1,\dots,m$. The integral in \eqref{boundary31} vanishes, unless $i=1$ and $j=t$, or $i=t$ and $j=1$. Without loss of generality we may assume that $j=1$.

If $t\leq m$, then integral is bounded by
\begin{align*}
\int\limits_{S_t(\eps,R)} |\log r_n|^2 \frac{d\rho_1 dr_1\cdots d\rho_{t-1} dr_{t-1}\cdot d\rho_t \cdot  d\rho_{t+1} dr_{t+1}\cdots  d\rho_p dr_p}{\psi^*(Q(Y)^p) \cdot r_{\alpha_1} \cdots r_{\alpha_\mu} },
\end{align*}
and $\alpha_1,\dots,\alpha_\mu\neq t$.
Since $\lambda_{\alpha_1},\dots,\lambda_{\alpha_\mu}$ are linearly independent, after possibly renaming indices, we may assume that $\lambda_t,\lambda_{\alpha_2},\dots,\lambda_{\alpha_\mu}$ are also linearly independent. Hence, it follows from Lemma \ref{cone2}a with $A=1-2/\mu$ (respectively directly from \eqref{impsi} if $\mu=1$) that
\[
\psi^*(Q(Y)^p)\gg |\log r_t|\cdot \left( |\log r_{\alpha_2}|\cdots |\log r_{\alpha_\mu}|\right)^{1/(\mu-1)}.
\]
Consequently, the integral in \eqref{boundary31} is bounded by
\begin{align*}
&\int\limits_{S_t(\eps,R)} |\log r_n|^2  \frac{d\rho_1 dr_1\cdots d\rho_{t-1} dr_{t-1}\cdot d\rho_t \cdot  d\rho_{t+1} dr_{t+1}\cdots  d\rho_p dr_p}{|\log r_t|^p\cdot \left( |\log r_{\alpha_2}|\cdots |\log r_{\alpha_\mu}|\right)^{p/(\mu-1)}\cdot r_{\alpha_1} \cdots r_{\alpha_\mu} } \stackrel{\eqref{intest3}}{\ll} |\log \eps|^{3-p}. 
\end{align*}
We obtain  \eqref{boundary31}, if $p>3$. In the case $p\leq 3$ our assumption on the Witt rank of $L$ implies that $K$ is anisotropic. We may apply Lemma \ref{cone2}b to deduce that the the integral in \eqref{boundary31} is actually $O(|\log \eps|^{-1})$ as in the proof of Lemma \ref{boundary1}.

If $t>m$, then $\alpha_1,\dots,\alpha_\mu\neq t$ and $\beta_1,\dots\beta_\nu\neq t$ in  \eqref{boundary31}.
Hence, by means of the estimate  $\psi^*(Q(Y))\gg 1$, we find that the integral in \eqref{boundary31} is bounded by
\[
\int\limits_{S_t(\eps,R)} |\log r_n|^2 \frac{d\rho_1 dr_1\cdots d\rho_{t-1} dr_{t-1}\cdot r_t d\rho_t \cdot  d\rho_{t+1} dr_{t+1}\cdots  d\rho_p dr_p}{r_{\alpha_1} \cdots r_{\alpha_\mu} }.
\]
In view of \eqref{intest3}, this is $\ll \eps |\log\eps|^C$ for some $C>0$, yielding \eqref{boundary31}.

We now consider the case that $m<\kappa\leq p$.
In the same way as above, one shows that
\[
\lim_{\eps \to 0} \int\limits_{S_t(\eps,R)} d\rho_\kappa\wedge \eta = 0,
\]
if $t\neq \kappa$. If $t=\kappa$, then
\[
\frac{1}{2\pi}\lim_{\eps \to 0} \int\limits_{S_\kappa(\eps,R)} d\rho_\kappa\wedge \eta 
= \int\limits_{q_\kappa=0} \eta.
\]
To obtain the absolute convergence of the integral on the right hand side we notice that only the $\widehat{dq_\kappa} \widehat{d\bar q_\kappa}$-component of $\psi^*(\Omega^{p-1})$ gives a contribution. Thus, by the argument of Lemma \ref{boundary1}, the integral 
is bounded by a sum of integrals of the form 
\begin{align*}
& \int\limits_{q_\kappa=0}\sigma \frac{|\log r_1|^2+\dots+|\log r_m|^2}{\psi^*(Q(Y)^p)} \frac{ \widehat{dq_\kappa} \widehat{d\bar q_\kappa}}{ r_{\alpha_1} \cdots r_{\alpha_\mu} r_{\beta_1}\cdots r_{\beta_{\nu}}}\\
&\ll \int\limits_{r_1=0}^R \cdots \int\limits_{r_{\kappa-1}=0}^R \cdot \int\limits_{r_{\kappa+1}=0}^R \cdots \int\limits_{r_{p}=0}^R  \frac{|\log r_1|^2+\dots+|\log r_m|^2}{\psi^*(Q(Y)^p)} \frac{dr_1\cdots dr_{\kappa-1}\cdot dr_{\kappa+1}\cdots dr_{p}}{ r_{\alpha_1} \cdots r_{\alpha_\mu}}\\
&\ll \int\limits_{r_1=0}^R \cdots \int\limits_{r_{m}=0}^R
\frac{|\log r_1|^2+\dots+|\log r_m|^2}{\psi^*(Q(Y)^p)} \frac{dr_1\cdots dr_{m}}{ r_{\alpha_1} \cdots r_{\alpha_\mu}}.
\end{align*}
The convergence of this integral follows in the same way as before.
\end{proof}

\begin{cor}
The analytic degree $\deg_{Y_\Gamma}(D)$ of a divisor $D$ on $Y_\Gamma$ only depends on the linear equivalence class of $D$.
\end{cor}
 

Theorem \ref{intro2} is very useful for determining the divisor of a given holomorphic modular form $F$ for $\Gamma$ of known weight. For instance, if we combine it with the explicit formula for the degrees of Heegner divisors on $Y_\Gamma$ (Proposition 4.8 in \cite{BrKue}), we find that the assumption that $\div(F)$ be a linear combination of Heegner divisors can be dropped from Theorem 13 in \cite{BK}.

\end{document}